\documentclass[12pt]{article}
\usepackage{amsbsy}
\usepackage{amssymb}
\usepackage{amsfonts}
\usepackage{amsmath,amscd}
\usepackage{amsopn}
\usepackage{amsthm}
\usepackage[applemac]{inputenc}
\usepackage[T1]{fontenc}

\newtheorem{theorem}{Theorem}
\newtheorem{corollary}{Corollary}

\newtheorem{lemma}{Lemma}
\theoremstyle{definition}
\newtheorem{remark}{Remark}

\theoremstyle{definition}
\newtheorem{example}{Example}
\newtheorem{definition}{Definition}
\newcommand{\p}{\mathbb{P}}

\newcommand{\e}{\mathbb{E}}

\title{Special, conjugate and complete
scale functions for spectrally negative L\'evy processes}
\author{A. E. Kyprianou\thanks{University of Bath, UK}\quad and  V. Rivero\thanks{Centro de Investigación en Matemáticas A.C., Mexico \& University of Bath, UK}}

\begin{document}

\maketitle
 \begin{abstract}
Following from recent developments in Hubalek and Kyprianou \cite{HK2007} the objective of this paper is to provide further methods for constructing new families of scale functions for spectrally negative L\'evy processes
which are completely explicit.  This is the result of an observation in the aforementioned paper which permits feeding the theory of  Bernstein functions directly into the Wiener-Hopf factorization for spectrally negative L\'evy processes. Many new, concrete examples of scale functions are offered although the methodology in principle delivers still more explicit examples than those listed.
\end{abstract}

\section{Spectrally negative L\'evy processes and scale functions}
Let $X=\{X_t : t\geq 0\}$ be a L\'evy process defined on a
filtered probability space $(\Omega , \mathcal{F}, \mathbb{P})$,
where $\{\mathcal{F}_t: t\geq 0\}$ is the filtration generated by $X$ satisfying the usual conditions.
For  $x\in \mathbb{R}$ denote by  $\mathbb{P}_{x}$ the law of $X$ when it is started at $x$ and
write simply  $\mathbb{P}_{0}=\mathbb{P}$. 
Accordingly we shall write $\mathbb{E}_x$ and $\mathbb{E}$ for the associated expectation operators. In this paper we shall
assume throughout that $X$ is {\it spectrally negative} meaning here that it has no positive jumps and that it is not the
negative of a subordinator. It is well known that the latter allows us to talk about the Laplace exponent $\psi(\theta):=
\log\mathbb{E}[e^{\theta X_1}]$ for $\theta\geq 0$ where in particular we have the L\'evy-Khintchine representation 

\begin{equation}
\psi(\theta) = -a\theta + \frac{1}{2}\sigma^2\theta^2 + \int_{(-\infty, 0)}(e^{\theta x} -1- x\theta\mathbf{1}_{\{x>-1\}})\Pi(d x)
\label{Laplace}
\end{equation}
where $a\in\mathbb{R}$, $\sigma\geq 0$ is the Gaussian coefficient and $\Pi$ is a measure concentrated on $(-\infty,0)$ satisfying $\int_{(-\infty,0)}(1\wedge x^2) \Pi(dx)<\infty$. 
The, so-called, L\'evy triple $(a,\sigma, \Pi)$ completely characterizes the process $X$.

For later reference we also introduce the 
function $\Phi: [0,\infty)\rightarrow [0,\infty)$ as 
the right inverse of $\psi$  on $(0,\infty)$ so that 
for all $q\geq 0$
\[
\Phi(q) =  \sup\{\theta\geq 0 : \psi(\theta) = q\}.
\]
Note that it is straightforward to show that  $\psi$ is a strictly convex
function which is zero at the origin and tends to infinity at infinity and hence
there are at most two solutions of the equation $\psi(\theta)=q$.

Suppose now we define the stopping times for each $x\in\mathbb{R}$
\[
\tau^+_x =\inf\{t> 0 : X_t >x\}\text{ and }\tau^-_x = \inf\{t> 0 : X_t < x\}.
\]

A fluctuation identity with a long history concerns the probability that $X$ exits an interval $[0,a]$ (where $a>0$) into
$(a,\infty)$ before exiting into $(-\infty,0)$ when issued at $x\in[0,a]$.  In particular it is known that 
\begin{equation}
\mathbb{E}_x(e^{-q\tau^+_a} \mathbf{1}_{\{\tau^+_a< \tau^-_0\}}) = \frac{W^{(q)}(x)}{W^{(q)} (a)}
\label{twosidedexit}
\end{equation}
where $x\in (-\infty,a]$, $q\geq 0$ and the function $W^{(q)} : \mathbb{R}\rightarrow [0,\infty)$ is defined up to a
multiplicative constant as follows. On $(-\infty,0)$ we have $W^{(q)}(x)=0$ and otherwise $W^{(q)}$ is a continuous function
(right continuous at the origin) with Laplace transform 
\[
\int_0^\infty e^{-\theta x}W^{(q)}(x)d x = \frac{1}{\psi(\theta) -q}\quad \text{ for }\theta>\Phi(q).
\label{LT}
\]

The functions $\{W^{(q)} : q\geq 0\}$ are known as {\it scale functions} and for convenience and consistency we write $W$ in place of $W^{(0)}$.
Identity (\ref{twosidedexit}) exemplifies the relation between scale functions for $q=0$ and the classical ruin problem. Indeed setting $q=0$ we have that $\mathbb{P}_x(\tau^-_0<\tau^+_a) = 1-W(x)/W(a)$ and so assuming that $\psi'(0+)>0$ and taking limits as $a\uparrow\infty$ it is possible to deduce that $W(\infty)^{-1} = \psi'(0+)$ and hence 

\[
\mathbb{P}_x(\tau^-_0<\infty)=1 - \psi'(0+)W(x).
\]
It is in this context of ruin theory that they 
make their earliest appearance in the works of \cite{Zol1964}, 
\cite{Tak1966} and then later either explicitly or implicitly in the work of \cite{Eme1973} 
and \cite{Kor1974,Kor1975}, \cite{Sup1976} and \cite{Rog1990}, 
the real value of scale functions as a  class with which one may express a whole range of fluctuation 
identities for spectrally negative L\'evy processes became apparent in the work of \cite{Chau94, Chau96} and  \cite{Ber1996, Ber1997} and an ensemble of subsequent articles; 
see for example \cite{Lam2000}, \cite{AKP2004}, 
\cite{P2003,Pis2004,Pis2005,Pis2006}, 
\cite{KP2006}, \cite{DK2006} and \cite{Don1991,Don2005,Don2007}. 
Moreover with the advent of these new fluctuation identities and a better understanding of the analytical 
properties of the function $W^{(q)}$ came the possibility of revisiting and solving a number of classical 
and modern problems from applied probability, but now with the underlying source of randomness being a 
general spectrally negative L\'evy processes. For example, in the theory of optimal 
stopping \cite{ACU2002, AKP2004} and \cite{Kyp2006}, in the theory of optimal 
control \cite{APP2007}, \cite{KP2007, RZ2007} and \cite{Ronnie2007}, in the theory of queuing and storage 
models \cite{DGM2004} and \cite{BBK2007}, in the theory of branching 
processes \cite{Bin1976}, \cite{huillet2003} and \cite{Lam2007}, in the theory of insurance risk and 
ruin \cite{CY2005}, \cite{KKM2004}, 
\cite{KK2006}, \cite{DK2006}, 
\cite{KP2006},  \cite{RZ2007} and \cite{Ronnie2007}, in the theory of credit risk 
\cite{HR2002} and \cite{KS2007} and in the theory of fragmentation \cite{Krell2007}. %See also \cite{Kyp2006} for a general overview of scale  functions and their applications.

Although scale functions are now firmly embedded within the theory of spectrally negative L\'evy 
processes and their applications, and although there is a reasonable understanding of how they 
behave analytically (see for example the summary in \cite{Kyp2006}), one of their main failings 
from a practical point of view until recently is that there are a limited number of concrete examples.  
None the less, Hubalek and Kyprianou \cite{HK2007} have been able to describe a parametric family of scale functions in explicit detail using a mathematical trick in which a spectrally negative L\'evy process is constructed having a particular pre-determined Wiener-Hopf factorization.

In this paper, we make use of the aforementioned trick and combine it with classical and recent developments in the potential analysis of subordinators. Stemming from the latter,  the theory of special Bernstein functions feed directly into the theory of scale functions and brings about the existence of pairs of spectrally negative L\'evy processes whose scale functions are conjugate to one another in an appropriate sense. 

Following an analytical exposition of the latter concepts, we are able to offer completely new explicit examples of conjugate pairs of scale functions.

According to Lemma 8.4 in \cite{Kyp2006}, for any $q>0,$ the $q$-scale function of a spectrally negative L\'evy process can be obtained via the $0$-scale function of a spectrally negative L\'evy process that drifts to $\infty,$ whose law is obtained by an exponential change of measure of the law of the original L\'evy process. From the theoretical point of view this allows us to restrict our study to  only the case $q=0.$ Nevertheless, from the point of view of applications this represents a restriction given that to make such construction tractable one needs to know explicitly the bivariate Laplace exponent associated to the descending ladder process, its potential measure and the inverse of the Laplace exponent $\psi$, which in general is an open problem. 

In the sequel we will assume, unless otherwise stated, that the spectrally negative L\'evy process does not drift to $-\infty.$ The reason for this is that the $0$-scale function of a L\'evy processes drifting to $-\infty$ can be deduced from the $0$-scale function corresponding to the associated L\'evy process conditioned to drift to $\infty.$ This will be made precise in Section \ref{genconst} below.    

\section{Descending ladder height and parent processes}
The principal idea for generating new examples of scale functions in Hubalek and Kyprianou \cite{HK2007}, which we borrow here, relies on constructing a spectrally negative L\'evy process around a given possibly killed subordinator  which plays the role of the descending ladder height process. For the convenience of the reader we devote a little time in this section reminding the reader of the meaning of a descending ladder height process and give the result of \cite{HK2007} in detail.

It is straightforward to show that the process $X-\underline{X}:=\{X_t-\underline{X}_t:t\geq 0\}$, 
where $\underline{X}_t := \inf_{s\leq t}X_s$, is a strong Markov process with state space $[0,\infty)$. 
Following standard theory of Markov local times (cf. Chapter IV of \cite{Ber1996}), it is possible to construct a  local time at zero 
for $X-\underline{X}$ which we henceforth refer to as $L=\{L_t : t\geq 0\}$. Its inverse process, $L^{-1}:=\{L^{-1}_t : t\geq 0\}$ where $L^{-1}_t =\inf\{s>0 : L_s >t\}$, 
is a (possibly killed) subordinator.
Sampling $X$ at $L^{-1}$ we recover the points of minima of $X$. If we define $H_t =X_{L^{-1}_t}$ 
when $L^{-1}_t<\infty,$ with $H_t  =\infty$ otherwise, then it is known that the 
process $H=\{H_t : t\geq 0\}$ is a (possibly killed) subordinator. The latter is known as 
the {\em descending ladder height process}. Moreover,  if $\Upsilon$ is the jump measure of $H$ then 
\[
\Upsilon(x,\infty) = e^{\Phi(0)x}\int_x^\infty e^{-\Phi(0)u}\Pi(-\infty,-u)d u\quad \text{for}\ x>0,
\]
see for example \cite{Vig2002}. 
Further, the subordinator has a drift component if and only if $\sigma>0$ in which case the 
drift is necessarily equal to $\sigma^2/2$. The killing rate of $H$ is given by the constant
$\mathbb{E}(X_1)\vee 0 = \psi'(0+)\vee 0$.%; see for example \cite[Exercise~6.5]{Kyp2006}. 
\ Observe that in the particular case where $\Phi(0)=0$ the jump measure of $H$ has a non-increasing density.

The starting point for the relationship between the descending ladder height process and scale functions 
is given by the Wiener-Hopf factorization. In `Laplace form', for spectrally negative L\'evy processes, this can be written as
\begin{equation}
\psi(\theta) = (\theta -\Phi(0))\phi(\theta),\qquad \theta\geq 0,
\label{WHfact}
\end{equation}
where $\phi(\theta) = -\log \mathbb{E}(e^{-\theta H_1})$.
In the special case that $\Phi(0)=0$, that is to say, the process $X$ does not drift to $-\infty$ or equivalently that $\psi'(0+)\geq 0$, it can be shown that the scale function $W$ describes the potential measure of $H$. In other words
\begin{equation}
\int_0^\infty d t \cdot \mathbb{P}(H_t \in d x) = W(d x)\qquad \text{ for }x\geq 0
\label{resolvent2}
\end{equation}
or equivalently
\begin{equation}
\int_0^\infty e^{-\theta x}W(dx) = \frac{1}{\phi(\theta)}\qquad \text{ for }\theta>0.
\label{equiv}
\end{equation}
We henceforth restrict the discussion to the case that $\Phi(0)=0,$ see however Section \ref{genconst}.

The following theorem, taken from Hubalek and Kyprianou \cite{HK2007}, shows how one may identify a spectrally negative L\'evy process $X$ (called the {\it parent process}) for a given descending ladder height process $H$. An equivalent result has been obtained in Proposition 7 in \cite{BRY04}. Note that the version of the theorem we present here constructs the parent process such that it does not drift to $-\infty$. However, this is not a necessary restriction in the original formulation of this result.
\begin{theorem}\label{inverse}
Suppose that $H$ is a subordinator, killed at rate $\kappa\geq 0$, with jump measure which is absolutely continuous with non-increasing density, say $\upsilon,$ and drift ${\rm d}$. 
%Suppose further that $\varphi\geq 0$ is 
%given such that $\varphi\kappa=0$. Then, assuming that 
%\[
%\int_1^\infty x\left|\frac{d^2 \Upsilon}{d x^2}(x)\right|d x <\infty
%\]
%in the case that $\varphi=0$,
Then there exists a spectrally negative L\'evy process  $X$ that does not drift to $-\infty$, henceforth referred to as the {\it parent process}, 
%such that for all $x\geq 0$, $\mathbb{P}(\tau^+_x <\infty)=e^{-\varphi x}$ 
whose descending ladder 
height process is  the process $H$. The L\'evy triple $(a,\sigma, \Pi)$ of the parent process 
is uniquely identified as follows.
The Gaussian coefficient is given by 
\[
\sigma = \sqrt{2{\rm d}}.
\]
 The L\'evy measure is given by 
\[
\Pi(-\infty, -x) = \upsilon(x),\quad \text{for}\ x>0.%\varphi\Upsilon(x,\infty)+-\frac{ d \Upsilon(x,\infty)}{d x} 
\label{density}
\]
Finally 
\[
a =\int_{(-\infty, -1)} x\Pi(d x)  - \kappa.
\]
%if $\varphi=0$ and otherwise when $\varphi>0$ we can establish the value of $a$ from the equation
%\[
%a\varphi = \frac{1}{2}\sigma^2\varphi^2 + \int_{(-\infty, 0)}(e^{\varphi x} -1- x\varphi\mathbf{1}_{\{x>-1\}})\Pi(d x).
%\]
The Laplace exponent of the parent process is also given by
\[
\psi(\theta) = \theta\phi(\theta)
\]
for $\theta \geq 0$ where $\phi(\theta) = - \log \mathbb{E}(e^{-\theta H_1})$. The parent process oscillates or drifts to $\infty$ according to whether $\phi(0)=0$ or $>0.$  
\end{theorem}
   
Note that when describing parent processes later on in this text, for practical reasons we shall prefer to specify the triple $(\sigma, \Pi, \psi)$ instead of $(a,\sigma, \Pi)$. However both triples provide an equivalent amount of information.

\section{Special and conjugate scale functions}\label{spec-conj-scl-fn}
In this section we introduce the notion of a special Bernstein functions and special subordinators and use the latter to justify the existence of pairs of so called  {\it conjugate scale functions} which have a particular analytical structure. We refer the reader to the lecture notes of Song and Vondra\v{c}ek \cite{SV2007} and the books of  Berg and Forst \cite{BCFG1975}  and Jacob \cite{J1} for a more complete account of the theory of Bernstein functions and their application in potential analysis.

Recall that the class of Bernstein functions coincides precisely with the class of Laplace exponents of possibly killed subordinators. That is to say, a general Bernstein function takes the form 
\begin{equation}
\phi(\theta) = \kappa + {\rm d}\theta + \int_{(0,\infty)} (1 - e^{-\theta x})\Upsilon (dx) \text{ for }\theta\geq 0
\label{Bernstein}
\end{equation}
where $\kappa\geq 0 $, ${\rm d}\geq 0$ and $\Upsilon$ is a measure concentrated on $(0,\infty)$ such that $\int_{(0,\infty)}(1\wedge x)\Upsilon(dx)<\infty$.

\begin{definition}
Suppose that $\phi(\theta)$ is a Bernstein function, then it is called a {\it special Bernstein function} if 
\begin{equation}
 \phi(\theta) = \frac{\theta}{\phi^*(\theta)},\qquad \theta\geq 0,
\label{SB}
\end{equation}
where $\phi^*(\theta)$ is another Bernstein function. Accordingly a possibly killed subordinator is called a special subordinator if its Laplace exponent is a special Bernstein function.
\end{definition}
Note that if $\phi$ is a special Bernstein function with representation as given in (\ref{SB}) then one says that the Bernstein function $\phi^*$ is conjugate to $\phi$. Moreover it is apparent from its definition that $\phi^*$ is a special Bernstein function and $\phi$ is conjugate to $\phi^*$.   In \cite{Haw1975} and \cite{SV2007a} it is shown that a sufficient condition for $\phi$ to be a special subordinator is that $\Upsilon(x,\infty)$ is log-convex on $(0,\infty)$.

For conjugate pairs of special Bernstein functions $\phi$ and $\phi^*$ we shall write in addition to (\ref{SB})
\begin{equation}
 \phi^*(\theta) = \kappa^* + {\rm d}^*\theta + \int_{(0,\infty)}(1 - e^{-\theta x})\Upsilon^*(dx),\qquad \theta\geq 0,
 \label{Bernstein*}
\end{equation}
where necessarily $\Upsilon^*$ is a measure concentrated on $(0,\infty)$ satisfying $\int_{(0,\infty)}(1\wedge x)\Upsilon^*(dx)<\infty$. One may express the triple $(\kappa^*, {\rm d}^*, \Upsilon^*)$ in terms of related quantities coming from the conjugate $\phi$. Indeed it is known that 
\[
 \kappa^* = \left\{
\begin{array}{ll} 
0 & \kappa>0 \\
\left({\rm d} + \int_{(0,\infty)} x\Upsilon(dx)\right)^{-1} & \kappa = 0
\end{array}
\right.
\]
and
\begin{equation}
 {\rm d}^* = \left\{\begin{array}{ll}
0 & {\rm d}>0 \text{ or } \Upsilon(0,\infty)=\infty \\
\left( \kappa + \Upsilon(0,\infty)\right)^{-1} & {\rm d}=0 \text{ and }\Upsilon(0,\infty)<\infty.
\end{array}\right.
\label{d*}
\end{equation}
 Which implies in particular that $\kappa\kappa^*=0={\rm dd}^*.$ In order to describe the measure $\Upsilon^*$ let us denote by 
 $W(dx)$ 
the potential measure of $\phi$.  (This choice of notation will of course prove to be no coincidence). Then we have that $W$ necessarily satisfies
\[
 W(dx) = {\rm d}^*\delta_0(dx) + \{\kappa^* + \Upsilon^*(x,\infty)\}dx \text{ for }x\geq 0.
\]
 Naturally, if $W^*$ is the potential measure of $\phi^*$ then we may equally describe it in terms of $(\kappa, {\rm d}, \Upsilon)$. A proof of these facts and other interesting results can be found in \cite{SV2006}.
Besides, it can be easily shown that a necessary and sufficient condition for a function to be a special Bernstein function is that its potential measure has a density on $(0,\infty)$
which is non-increasing and integrable in the neighborhood of the origin; see e.g. \cite{Ber1997b} Corollaries 1 and 2 for a proof of this fact.

\bigskip

We are interested in constructing a parent process whose descending ladder height process is a special subordinator. The first parts of the following theorem and corollary are now evident given the discussion in the current and previous sections.

\begin{theorem}\label{specialscale}
 For conjugate special Bernstein functions $\phi$ and $\phi^*$ satisfying (\ref{Bernstein}) and (\ref{Bernstein*}) respectively where  $\Upsilon$ is absolutely continuous with non-increasing density, there exists a spectrally negative L\'evy process that does not drift to $-\infty$, 
whose Laplace exponent is described by 
\begin{equation}
 \psi(\theta) = \frac{\theta^2}{\phi^*(\theta)} = \theta\phi(\theta)\text{ for }\theta\geq 0
 \label{parent}
\end{equation}
and
whose scale function is a concave function and is given by 
\begin{equation}
 W(x) = {\rm d}^* + \kappa^*x + \int_0^x \Upsilon^*(y,\infty)dy.
 \label{specialW}
\end{equation}
Conversely, if $\psi$ is the Laplace exponent of a spectrally negative L\'evy process that does not drift to $-\infty$ and its associated scale function, $W,$ is a concave function then there exists a pair of conjugate Bernstein functions $\phi$ and $\phi^*,$ satisfying (\ref{Bernstein}) and (\ref{Bernstein*}) respectively, where  $\Upsilon$ is absolutely continuous with non-increasing density, and such that the relations (\ref{parent}) and (\ref{specialW}) hold.
\end{theorem}

\begin{proof}
Only the second part of the theorem needs a proof. Let $\psi$ and $W$ be as described in the second part of the statement of the theorem and let $\phi$ denote the Laplace exponent of the descending ladder height subordinator associated to the spectrally negative L\'evy process with Laplace exponent $\psi$. The latter and former functions are related via (\ref {WHfact})
 with $\Phi(0)=0$. By an integration by parts it follows that    
\begin{equation*}
\begin{split}
\theta\int^\infty_{0}e^{-\theta x}W(x)d x=W(0+)+\int^\infty_{0}W'(y)e^{-\theta y}d y=\frac{1}{\phi(\theta)},\qquad \theta\geq 0,  
\end{split}
\end{equation*} where $W'$ denotes the first derivative of $W,$ which exists almost everywhere because of the concavity of $W.$
This implies in particular that in $(0,\infty)$ the potential measure of the descending ladder height subordinator has a  density which is non-increasing as well as an atom at zero of size $W(0+)$. Furthermore, again by an integration by parts it follows that 
\begin{equation}\label{scalebernstein}
\begin{split}
&\theta^2\int^\infty_{0}e^{-\theta x}W(x)d x\\
&=\theta W(0)+\theta\int^\infty_{0}W'(y)e^{-\theta y}d y\\
&=W'(\infty)+\theta W(0)+\int^\infty_{0}\left(1-e^{-\theta y}\right)d \left(-W'(y)\right),\qquad \theta \geq 0. 
\end{split}
\end{equation} 
So, that the function $\phi^*$ defined by $$\phi^*(\theta)=\frac{\theta^2}{\psi(\theta)}=\frac{\theta}{\phi(\theta)},\qquad \theta\geq 0,$$ is a special Bernstein function conjugated to $\phi.$
\end{proof}
 Note that the proof of the converse statement in Theorem \ref{specialscale} says a little more than is claimed. Indeed we have that the potential measure of the subordinator with Laplace exponent $\phi^*,$ in $(0,\infty),$ admits a density which is decreasing and convex. To see this, one should recall that the potential measure associated to $\phi^*$ has a decreasing density which is given by the tail L\'evy measure of $\phi,$ and since $\psi(\theta)=\theta\phi(\theta),$ $\theta\geq0$ it follows by an integration by parts that the L\'evy measure of $\phi$ has a decreasing density.

\bigskip

The assumptions of the previous theorem require only that  the L\'evy and potential measures associated to $\phi$ have a non-increasing density in $(0,\infty)$, respectively; this condition on the potential measure is equivalent to the existence of $\phi^*.$ If in addition it is  assumed that the potential density be a convex function, in light of the representation (\ref{specialW}), we can interchange the roles of $\phi$ and $\phi^*,$ respectively, in the previous Theorem. The key issue to this additional assumption is that it ensures the absolute continuity of $\Upsilon^*$ with a non-increasing  density. We thus  have the following Corollary.    

\begin{corollary}\label{corrspecialscale}
If conjugate special Bernstein functions $\phi$ and $\phi^*$ exist satisfying (\ref{Bernstein}) and (\ref{Bernstein*}) such that both $\Upsilon$ and $\Upsilon^*$ are absolutely continuous with non-increasing densities, then there exist a pair of scale functions $W$ and $W^*$, such that $W$ is concave, its first derivative is a convex function,  (\ref{specialW}) is satisfied, and 
\begin{equation}
 W^*(x) = {\rm d} + \kappa x + \int_0^x \Upsilon(y,\infty)dy
 \label{specialW*}
\end{equation}
whose parent processes are given by (\ref{parent}) and
\begin{equation}\label{parent*}
\psi^*(\theta) = \frac{\theta^2}{\phi(\theta)} = \theta\phi^*(\theta).
\end{equation}
Conversely, if $\psi$ is the Laplace exponent of a spectrally negative L\'evy process such that its associated scale function $W$ is a concave function whose first derivative is a convex function, then there exists a pair of conjugate Bernstein functions $\phi$ and $\phi^*$ and a function $\psi^*,$ such that $\psi^*$ is the Laplace exponent of a spectrally negative L\'evy process, and $\psi$ (respectively $\psi^*$) is related to $\phi^*$ (respectively to $\phi$) by equation (\ref{parent}) (respectively by (\ref{parent*})) and the scale functions associated to $\psi$ and $\psi^*$ satisfy equations (\ref{specialW}) and (\ref{specialW*}), respectively.  
\end{corollary}

\begin{proof}
The first part of the proof is a simple consequence of Theorem \ref{specialscale}. The converse part follows from the calculations used in the proof of Theorem \ref{specialscale}. Indeed, if $\psi$ and $W$ are as described in the converse statement of the corollary then it follows from Theorem \ref{specialscale} that there exists a pair of conjugate Bernstein functions $\phi$ and $\phi^*$ and by hypotheses and equation (\ref{scalebernstein}) that the respective L\'evy measure of $\phi$ and $\phi^*$ have non-increasing densities. That the function $\psi^*$ defined by $\psi^*(\theta)=\theta\phi^*(\theta),$ $\theta\geq 0,$ is the Laplace exponent of a spectrally negative L\'evy process is then a consequence of Theorem \ref{inverse}. It follows from Theorem \ref{specialscale} and the uniqueness of Laplace transform that the scale functions associated to $\psi$ and $\psi^*$ have the claimed properties.\end{proof}

There are a number of remarks that are worth making regarding the above theorem in the setting of conjugate special Bernstein functions.

\begin{enumerate}
\item For obvious reasons we shall henceforth refer to the scale functions identified in (\ref{Bernstein}) and (\ref{Bernstein*}) as {\it special scale functions}.

\item Similarly, when $W$ and $W^*$ exist then we refer to them as {\it conjugate (special) scale functions} and their respective parent processes are called {\it conjugate parent processes}. This conjugation can be seen by noting that thanks to (\ref{SB}).
\[
W*W^*(dx) = dx.
\]

\item It is known (cf. Chapter 8 of \cite{Kyp2006}) that any scale function has a discontinuity at the origin if and only if the parent process has paths of bounded variation.  That is to say, in the representation (\ref{Laplace}), $\sigma=0$ and $\int_{(-1,0)}|x|\Pi(dx)<\infty$. Taking account of the description in Theorem \ref{inverse} one may easily deduce that  in terms of the descending ladder height process this is equivalent with the fact that $\Upsilon(0,\infty)<\infty$ and ${\rm d}=0$ which, within the context of Theorem \ref{specialscale},  is equivalent to the case that ${\rm d}^*>0$ in (\ref{d*})%\ref{Bernstein}) 
\ as predicted by the general theory.

\item Another known general property of scale functions is that $W'(0+)<\infty$ if and only if, in the representation (\ref{Laplace}), $\sigma>0$ or $\sigma=0$ {\it and} $\Pi(-\infty,0)<\infty$. See e.g. Exercise 8.5 in \cite{Kyp2006}. The latter conditions in terms of the descending ladder height process are respectively equivalent to ${\rm d}>0$ or ${\rm d}=0$ {\it and} $d\Upsilon(0+)/dx<\infty.$ 
%In the current context of Bernstein functions one may also see why this is the case. It is a general fact that if ${\rm d}>0$ then, since $W$ is a potential measure of a subordinator with drift, $W'(0+)=1/{\rm d}<\infty$ (cf. Theorem III.5 in \cite{Ber1996}). Moreover under the circumstances at hand, if $d\Upsilon(0+)/dx<\infty$ {\it and } ${\rm d}=0$ then necessarily $H$ is a pure jump compound Poisson (possibly killed) subordinator. Since $1/\phi(\theta)$ corresponds to the Laplace transform of $W(dx)$, it follows by taking $\theta\uparrow\infty$ that $W(0+)+W'(0+) = 1/(\kappa+\Upsilon(0,\infty))<\infty$. Conversely suppose that $W'(0+)<\infty$. Then from the form of the scale function given in Theorem \ref{specialscale} it follows that $\kappa^*+\Upsilon^*(0,\infty)<\infty$ which in turn, by (\ref{d*}), implies that ${\rm d}>0$.

\end{enumerate}

\section{Tilting and parent processes drifting to $-\infty$}\label{genconst}
In this section we present two methods for which, given a scale function and associated parent process, it is possible to construct further examples of scale functions by appealing to two procedures. 

The first method relies on the following facts concerning translating the argument of a given Bernstein function. 
% constitutes a translation of a given Bernstein function in its argument which can be seen as  a form of exponential tilting. 
%The following Lemma  together with the previous results will prove to be a useful tool to construct further explicit examples of scale functions from established examples. Its proof is a simple consequence from the theory of exponential tilting for L\'evy processes.
\begin{lemma}\label{exptilting0}
Let $\phi$ be a special Bernstein function with representation given by (\ref{Bernstein}). Then for any $\beta\geq 0$ the function $\phi_{\beta}(\theta)=\phi(\theta+\beta),$ $\theta\geq 0,$ is also a special Bernstein function with killing term $\kappa_{\beta}=\phi(\beta),$ drift term ${\rm d}_\beta = {\rm d}$ and L\'evy measure $\Upsilon_\beta(dx)=e^{-\beta x}\Upsilon(d x),$ $x>0.$ Its associated potential measure, $W_\beta$, has a decreasing density  in $(0,\infty)$ such that  $W_{\beta}(dx)=e^{-\beta x}W'(x)dx,$ $x>0,$ where $W'$ denotes the density of the potential  measure associated to $\phi.$ Moreover, let $\phi^*$ and $\phi^{*}_{\beta},$ denote the conjugate Bernstein functions of $\phi$ and $\phi_{\beta},$ respectively. Then the following identity
\begin{equation}\label{bernsteintilting}
\phi^*_{\beta}(\theta)=\phi^*(\theta+\beta)-\phi^*(\beta)+\beta\int^\infty_{0}\left(1-e^{-\theta x}\right)e^{-\beta x}W'(x)d x,\qquad \theta\geq 0,
\end{equation}
holds.
\end{lemma}
\begin{proof} We will use the theory of change of mesure for L\'evy processes; for background on this topic see e.g. \cite{Sat1999} Section 33. Assume that under $\p,$ $H$ is a subordinator with Laplace exponent $\phi.$  It follows from the independence and homogeneity of the increments of $H$ that the process $M_{t}=e^{-\beta H_{t}},$ $t\geq0,$  is a submartingale in the natural filtration generated by $H$. Moreover, $\lim_{t\to 0+}\e\left(e^{-\beta H_{t}}\right)=1.$  It follows that the function $h_{\beta}(x):=e^{-\beta x}, x\in[0,\infty],$ under the assumption that $e^{-\infty}=0,$ is an excessive function  for the semigroup of the process $H.$ We denote by $\p^{(\beta)}$ the law of the $h$-transform of $\p$ via the excessive function $h_{\beta},$ that is $\p^{(\beta)}$ is the unique measure over the space of right continuous left-limited with lifetime paths such that $$\p^{(\beta)}=M_{t}\p\qquad \text{over}\ \sigma(H_{s}, s\leq t),\ \text{for}\ t\geq 0.$$ It is easily verified that under $\p^{(\beta)}$ the law of $H$ is that of a subordinator killed at rate $\phi(\beta).$ Indeed, as $\p^{(\beta)}$ is locally absolutely continuous with respect to $\p$ it follows that under $\p^{(\beta)},$ $H$ has non-decreasing paths and for $s,t\geq 0$
\begin{equation*}
\begin{split}
&\e^{(\beta)}\left(e^{-\theta(H_{t+s}-H_{t})}1_{A}\right)\\
&=\e\left(e^{-(\theta+\beta)(H_{t+s}-H_{t})}e^{-\beta H_{t}}1_{A}\right)\\
&=\e\left(e^{-(\theta+\beta)H_{s}}\right)\e(1_{A}e^{-\beta H_{t}})\\
&=e^{-s\phi(\theta+\beta)}\e^{(\beta)}(1_{A}),
\end{split}
\end{equation*}  
for every $\theta\geq 0$ and $A\in \sigma(H_{u}, u\leq t).$ Which proves at once that under $\p^{(\beta)},$ $H$ has independent and homogeneous increments, its Laplace exponent is given by $\phi_{\beta}(\cdot)=\phi(\beta+\cdot),$ and its killing term is thus $\phi_{\beta}(0)=\phi(\beta).$ Moreover, the equality 
\begin{equation*}
\begin{split}
\phi_{\beta}(\theta)&=\phi(\beta)+\theta {\rm d}+\int^\infty_{0}(e^{-\beta x}-e^{-(\beta+\theta)x})\Upsilon(dx)\\&=\phi(\beta)+\theta {\rm d}+\int^\infty_{0}(1-e^{-\theta x})e^{-\beta x}\Upsilon(dx),\qquad \theta\geq 0,
\end{split}
\end{equation*}
justifies the description of the characteristics of $\phi_{\beta}$ claimed in the Lemma~\ref{exptilting0}. Furthermore, the potential of the process $H$ under $\p^{(\beta)}$ is given by 
\begin{equation*}
\begin{split}
&\e^{(\beta)}\left(\int^\infty_{0}1_{\{H_{t}\in dx\}}dt\right)=\int^\infty_{0}\e^{(\beta)}\left(1_{\{H_{t}\in dx\}}\right)dt
=\int^\infty_{0}\e\left(1_{\{H_{t}\in dx\}}e^{-\beta H_{t}}\right)dt\\
&=\e\left(\int^\infty_{0}e^{-\beta H_{t}}1_{\{H_{t}\in dx\}}dt\right)={\rm d^*}\delta_{0}(dx)+e^{-\beta x}W'(x)dx, 
\end{split}
\end{equation*}
for $x\geq 0,$ owing to Fubini's Theorem. Given that the function $x\mapsto e^{-\beta x}W'(x),$ $x>0$ is a decreasing function, it follows that $\phi_{\beta}$ is a special Bernstein function. Furthermore, following the description in Section~\ref{spec-conj-scl-fn} the characteristics of its conjugate, $\phi^*_{\beta}$ are given as follows: its killing term is $\kappa^*_{\beta}=0,$ as $k_{\beta}=\phi(\beta)>0,$ the tail of its L\'evy measure is given by $$\Upsilon^{*}_{\beta}(x,\infty)=e^{-\beta x}W'(x),\qquad x>0,$$ and its drift equals ${\rm d}^*$ as $${\rm d}^*_{\beta}=\lim_{\theta\to\infty}\frac{\phi^*_{\beta}(\theta)}{\theta}=\lim_{\theta\to\infty}\frac{1}{\phi(\beta+\theta)}=\lim_{\theta\to\infty}\frac{1}{\phi(\theta)}={\rm d}^*,$$ in the obvius notation. The description in~(\ref{bernsteintilting}) follows by bare-hands calculations using the latter facts.\end{proof}
Note in particular that if $\Upsilon$ has a non-increasing density then so does $\Upsilon_\beta$. Moreover, if $W'$  convex (equivalently $\Upsilon^*$ has a non-increasing density) then $W'_\beta$ is  convex (equivalently $\Upsilon^*_\beta$ has a non-increasing density). These facts lead us to the following Lemma.
\begin{lemma}\label{exptilting}
If conjugate special Bernstein functions $\phi$ and $\phi^*$ exist satisfying (\ref{Bernstein}) and (\ref{Bernstein*}) such that both $\Upsilon$ and $\Upsilon^*$ are absolutely continuous with non-increasing densities, then there exist conjugate parent processes with Laplace exponents 
\begin{equation*}
 \psi_{\beta}(\theta)=\theta\phi_{\beta}(\theta)\text{ and }  \psi_{\beta}^*(\theta)=\theta\phi^*_{\beta}(\theta),\qquad \theta\geq 0.
\end{equation*}
whose respective scale functions are given by 
$$
W_{\beta}(x)={\rm d}^*+\int^x_{0}e^{-\beta y}\Upsilon^*(y,\infty)d y,\qquad x\geq 0.
$$
and
\begin{equation}\label{conjscaletilted}
W^*_{\beta}(x)={\rm d}+\phi(\beta)x+\int^x_{0}\left(\int^\infty_{y}e^{-\beta z}\Upsilon(d z)\right)d y.
\end{equation} 
using obvious notation. 
\end{lemma}

\bigskip

The second method builds on the latter to construct examples of scale functions whose parent processes drift to $-\infty.$ %It is worth observing that the process we will next construct can be seen as an auxiliary parent process constructed as in Section \ref{inverse}, which drifts to $\infty$, conditioned to drift to $-\infty$.

Suppose that  $\phi$ is a Bernstein function such that $\phi(0)=0,$ its associated L\'evy measure has a decreasing density and let $\beta>0$. Theorem \ref{inverse}, as stated in its more general form in \cite{HK2007}, says that there exists a parent process, say $X,$ that drifts to $-\infty$ such that its Laplace exponent $\psi$ can be factorized as $$\psi(\theta)=(\theta-\beta)\phi(\theta),\qquad \theta\geq 0.$$ It follows that $\psi$ is a convex function and $\psi(0)=0=\psi(\beta),$ so that $\beta$ is the largest positive solution to the equation $\psi(\theta)=0.$ Now, let $W_{\beta}$ be the $0$-scale function of the spectrally negative L\'evy process, say  $X_{\beta},$ with Laplace exponent 
$\psi_{\beta}(\theta):=\psi(\theta+\beta),$ for $\theta\geq 0.$ It is known that the L\'evy process $X_{\beta}$ is obtained by an exponential change of measure and can be seen as the L\'evy process $X$ conditioned to drift to $\infty,$ see chapter VII in \cite{Ber1996}. Thus the Laplace exponent $\psi_{\beta}$ can be factorized as $\psi_{\beta}(\theta)=\theta\phi_{\beta}(\theta),$ for $\theta\geq 0,$ where, as before, $\phi_{\beta}(\cdot):=\phi(\beta+\cdot).$ 
%Observe that $\phi_{\beta}$ is again a Bernstein function such that $\kappa_{\beta}=\phi(\beta),$ $d_{\beta}=d$ and $$\frac{\Upsilon_{\beta}(dx)}{dx}=e^{-\beta x}\frac{\Upsilon(dx)}{dx},\qquad x>0,$$ so that the L\'evy measure of $\phi_{\beta}$ has a decreasing density, in the obvious notation. 
It follows from Lemma 8.4 in \cite{Kyp2006}, that the $0$-scale function of the process with Laplace exponent $\psi$ is related to $W_{\beta}$ by $$W(x)=e^{\beta x}W_{\beta}(x),\qquad x\geq 0.$$ 
The above considerations thus lead to the following result which allows for the construction of a second parent process and associated scale function over and above the pair described in Theorem \ref{specialscale}.

\begin{lemma}\label{lemdriftneg} Suppose that $\phi$ is a special Bernstein function satisfying (\ref{Bernstein}) such that $\Upsilon$ is absolutely continuous with non-increasing density and $\kappa=0$. Fix $\beta>0$. Then there exists a parent process with Laplace exponent 
\[
\psi(\theta) = (\theta-\beta)\phi(\theta), \quad \theta\geq 0
\]
whose associated scale function is given by 
$$W(x)={\rm d}^*e^{\beta x}+e^{\beta x}\int^x_{0}e^{-\beta y}\Upsilon^*(y,\infty)dy,\qquad x\geq 0,$$ where we have used our usual notation. 
\end{lemma}
Now, when we assume furthermore that the potential density associated to $\phi$ is a decreasing and convex function, or equivalently that the L\'evy measure of $\phi^*$ has a decreasing density, there are three choices for a conjugate parent process. The first, is the one appearing in Lemma \ref{exptilting} with Laplace exponent given by $\psi^*_{\beta}(\theta)=\theta\phi^*_{\beta}(\theta),$ for $\theta\geq 0,$ and its scale function, $W^*_{\beta},$ is described in equation (\ref{conjscaletilted}). This parent process drifts to $\infty.$ The scale functions $W_{\beta}$ and $W^*_{\beta}$ are conjugated in the sense described in Remark 2 in Section \ref{spec-conj-scl-fn}, which implies that $$d(e^{-\beta x}W(x))*dW^*_{\beta}(x)=dx, \qquad x\geq 0.$$ The second and third candidate parent process, based on $\phi^*(\theta)$ are the one drifting to $-\infty,$ constructed using the formulation above, and the one which drifts to $\infty$ described in Corollary \ref{corrspecialscale}, accordingly as $\phi^*(0)=0$ or $\phi^*(0)>0,$ respectively. For these parent processes the respective scale functions are such that 
$$e^{-\beta x}dx=\begin{cases}
d\left(e^{-\beta x}W(x)\right)*d\left(e^{-\beta x}W^*(x)\right),& \text{if}\ \phi^*(0)=0,\\
d(e^{-\beta x}W(x))*(e^{-\beta x}W^*(dx)),& \text{if}\ \phi^*(0)>0,\end{cases} \qquad x\geq 0.$$
\begin{remark}
Observe that the construction explained in this section can be performed as soon as there exists a $\beta>0$ such that the function $\theta\mapsto \phi(\beta+\theta)$ is a special Bernstein function with a non-increasing L\'evy density. Which in view of the calculations carried in the proof of Lemma \ref{exptilting0} could occur without $\phi$ being a special Bernstein function in itself.
\end{remark}

\section{Complete scale functions}\label{completeclass}

We begin by introducing the notion of a complete Bernstein function with a view to constructing scale functions whose parent processes are derived from descending ladder height processes with Laplace exponents which belong to the class of complete Bernstein functions.

\begin{definition}\rm
A function $\phi$ is called {\it complete Bernstein function} if there exists an auxiliary Bernstein function $\eta$ such that 
\begin{equation}
\phi(\theta) = \theta^2 \int_{(0,\infty)} e^{-\theta x}\eta(x)dx.
\label{cB}
\end{equation}
\end{definition}
It is well known that a complete Bernstein function is necessarily a special Bernstein function (cf. \cite{J1}) and in addition, its conjugate is also a complete Bernstein function. Moreover, from the same reference one finds that a ne\-cessary and sufficient condition for $\phi$ to be complete Bernstein is that $\Upsilon$ satisfies for $x>0$
\[
\Upsilon(dx) = \left\{\int_{(0,\infty)} e^{- xy}\gamma(dy)\right\} dx
\]
where 
$
\int_{(0,1)}\frac{1}{y}\gamma(dy) + \int_{(1,\infty)}\frac{1}{y^2}\gamma(dy) <\infty.
$
Equivalently $\Upsilon$ has a completely monotone density.
Another necessary and sufficient condition is that the potential measure associated to $\phi$ has a density on $(0,\infty)$ which is completely monotone, this is a result due to Kingman \cite{Kin1967} and Hawkes \cite{Haw1976}. The class of infinitely divisible laws and subordinators related to this type of Bernstein functions has been extensively studied by several authors, see e.g. \cite{bondesson}, \cite{thorin}, \cite{rosinski}, \cite{donatiyor}, \cite{lancelotroynetteyor} and the references therein. 

Since necessarily $\Upsilon$ is absolutely continuous with  a completely monotone density, it follows that any subordinator whose Laplace exponent is a complete Bernstein function may be used in conjunction with Corollary \ref{corrspecialscale}. The following result is now a straightforward application of the latter and the fact that from  (\ref{cB}), any Bernstein function $\eta$ has a Laplace transform $(\theta^2/\phi(\theta))^{-1}$ where $\phi$ is complete Bernstein.

\begin{corollary}\label{completescale} Let $\eta$ be any Bernstein function and suppose that $\phi$ is the complete Bernstein function associated with the latter via the relation (\ref{cB}). Write $\phi^*$ for the conjugate of $\phi$ and  $\eta^*$ for the Bernstein function associated with $\phi^*$ via the related (\ref{cB}). Then 
\[
W(x) = \eta^*(x) \text{ and } W^*(x) = \eta(x),\qquad x\geq 0.
\]
are conjugate  scale functions with conjugate parent processes whose Laplace exponents are given by 
\[
\psi(\theta) = \frac{\theta^2}{\phi^*(\theta)}  = \theta\phi(\theta) \text{ and }\psi^*(\theta) = \frac{\theta^2}{\phi(\theta)} = \theta\phi^*(\theta),\qquad \theta\geq 0.
\]
\end{corollary}
%Note that this result is very closely related to a classical result which states that a Bernstein function has L\'evy density which is completely monotone if and only if its associated potential density is completely monotone. See \cite{Kin1967, Haw1976}.

We conclude this section with some remarks about the above corollary.
\begin{enumerate}

\item For notational consistency we call the pair $W$ and $W^*$ {\it (conjugate) complete scale functions} and their respective parent processes {\it (conjugate) complete parent processes}

\item In essence only part of the above corollary is of practical value. That is to say, any given Bernstein function $\eta$ is a scale function whose parent process is the spectrally negative L\'evy process whose Laplace exponent is given by $\psi^*(\theta) = \theta^2/\phi(\theta)$ where $\phi$ is given by (\ref{cB}).

\item Observe that any given completely monotone function say $p:[0,\infty[\to [0,\infty],$ such that $\int^1_{0}p(s)d s<\infty,$ may be seen as the potential density of a subordinator. This is due to the fact that any such function is log-convex, owing to H\"older's inequality, and then by a result due to Hirsch \cite{hirsch} there exists a subordinator, say $H,$ whose potential measure admits $p$ as a  density in $(0,\infty),$ see also \cite{SV2007a} for a recent proof of the latter fact. Let $\phi$ be the Laplace exponent of $H,$ this is such that the tail of its L\'evy measure is a completely monotone function that we will denote by $p^*.$ Thus, $\int^1_{0}p^*(s)d s<\infty$ and 
 $$
 \frac{\phi(\theta)}{\theta}={\rm d}+\frac{\kappa}{\theta}+\int^{\infty}_{0}e^{-\theta y}p^*(y)d y,\qquad \theta\geq 0.$$
 It follows from Corollary \ref{completescale} that there exists a spectrally negative L\'evy process with Laplace exponent $\psi(\theta)=\theta\phi(\theta)$ for $\theta\geq 0,$ its associated scale function is a Bernstein function, and can be represented as
$$
W(x)={\rm d}^*+\kappa^*x+\int^\infty_{0}(1-e^{-xy})\frac{\gamma(d y)}{y},\qquad x\geq 0;
$$ 
under the assumption that $p$ admits the representation 
$$
p(x)=\kappa^*+\int^\infty_{0}e^{-xy}\gamma(d y), \qquad x\geq 0,
$$ 
where $\kappa^*\geq 0$ and $\gamma$ is a measure over $(0,\infty)$ such that $$\int^1_{0}\gamma(d y)+\int^\infty_{1}\frac{\gamma(d y)}{y}<\infty,\quad \text{equivalently} \int^1_{0}p(t)d t<\infty.$$  For the respective conjugates we have that
 $$
 \phi^*(\theta)=\kappa^*+{\rm d}^*\theta+\int^\infty_{0}\frac{\theta}{\theta+y}\gamma(d y),\quad \theta\geq 0,
 $$ 
$$\sigma^*=\sqrt{2{\rm d}^*},\quad a^*= \int^\infty_{0}ye^{-y}\gamma(d y)+p(1)-\kappa^*,$$ 
$$\Pi^*(-\infty,-x)=\int^\infty_{0}ye^{-xy}\gamma(d y),\qquad x>0,$$ and therefore, 
$$W^*(x)={\rm d}+\kappa x+\int^x_{0}p^*(y)d y,\qquad x\geq 0.$$ 
Finally, owing to Remark 2 in Section \ref{spec-conj-scl-fn}, if ${\rm d}=0={\rm d}^*,$ the functions $p$ and $p^*$ are related by the Volterra-type equation 
$$\int^y_{0}p(x)p^*(y-x)d x=1=\int^y_{0}p^*(x)p(y-x)d x,\qquad y\geq 0.$$
\item Another source of examples comes from the observation that if $\varphi$ is a complete Bernstein function then the function $\phi(\theta)=\left(\varphi(1/\theta)\right)^{-1},$ for $\theta\geq 0,$ also is a Bernstein function. This assertion is easily proved using that the potential density associated to $\varphi$ is a completely monotone function. Thus, given a pair of complete conjugate Bernstein functions $\varphi,\varphi^*$ the functions $\phi,$ as defined above, and $\phi^*$ constructed analogously, form also a conjugate pair of complete Bernstein functions. So that having knowledge of the parent processes and scale functions associated to $\varphi$ and $\varphi^*,$ respectively, one can construct a new family of conjugate parent processes and scale functions.    
\end{enumerate}

\section{Concrete examples}
The previous sections have essentially consisted of  re-dressing the theory of Bernstein functions in the language of scale functions. In this section we justify the value of the previous exposition by offering a large cache of remarkably explicit examples.

\subsection*{Example 1}
Our first example will not be as interesting as other examples and has been included principally for the purpose of illustrating how the theory works in the context of an `old favorite'.

Consider the, apparently trivial, Bernstein function
\[
\phi(\theta) = \kappa + {\rm d}\theta
\]
where ${\rm d}, \kappa>0$. This is the Laplace exponent of the ladder height process with a parent process consisting of a Brownian motion with coefficient $\sqrt{2{\rm d}}$ drift at rate $\kappa$,
\[
\psi(\theta) = \kappa\theta + {\rm d}\theta^2,\qquad \theta\geq 0.
\]
This process is a diffusion and gives us an example where its scale function as a spectrally negative L\'evy process coincides precisely with its scale function as a diffusion. We have by a simple Laplace inversion
\[
W(x) = \frac{1}{\kappa} (1 - e^{- x\kappa/{\rm d }  }), \qquad x\geq0.
\]
This tells us that ${\rm d}^*=0$, $W'(x)=\kappa^*+\Upsilon^*(x,\infty) = {\rm d}^{-1}e^{-x\kappa/{\rm d}}$ and hence $\kappa^*=0$. Regarding the conjugate parent process, it is clear that 
\[
\Pi^*(-\infty, -x) = \frac{\kappa}{{\rm d}^2}e^{-x\kappa/{\rm d}}, \qquad x>0.
\]
and so the latter has a compound Poisson jump structure with negative exponentially distributed jumps having parameter $\kappa/{\rm d}$ and arrival rate $\kappa/{\rm d}^2$.
Since $\kappa^*={\rm d}^*=0$ the conjugate parent process is an oscillating process with no Gaussian component and hence one may write down from this information directly
\[
\psi^*(\theta) = \frac{1}{\rm d}\theta - \frac{\kappa}{{\rm d}^2}\left(1  - \frac{\kappa/{\rm d}}{\theta + \kappa/{\rm d}}\right).
\]
After a little algebra one finds that this coincides with the expected expression given by $\psi^*(\theta) = \theta^2/\phi(\theta)$. It follows from Theorem \ref{specialscale} that the scale function associated to $\psi^*$ is given by $$W^*(x)={\rm d}+\kappa x,\qquad x\geq 0.$$ So that $W^*$ is an ultimately linear scale function, thus its associated potential density is ultimately constant and so it is closely related to the potential measures appearing in Section 3 in \cite{SV2007a}.  Finally by the Continuity Theorem for Laplace transforms it follows that the the scale functions are continuous in the L\'evy triple $(\kappa,{\rm d},\Upsilon).$ For this reason as $\kappa\to 0$ the scale functions $W$ and $W^*$ converge towards the scale functions $x\mapsto x/{\rm d},$ and $x\mapsto {\rm d},$ for $x>0,$ respectively. Which provide two more examples of ultimately linear scale functions.   
%\end{example}

\subsection*{Example 2}
%\begin{example} 
Let $\beta,c>0$, $\nu\geq 0$ and  $\lambda\in(0,1)$. We claim that % and $\xi$ be the subordinator whose Laplace exponent $\phi\re^+\to\re^+,$ is given by the formula, 
 \[
 \phi(\theta)=\frac{c\beta\theta\Gamma(\nu+\beta\theta)}{\Gamma(\nu+\beta\theta+\lambda)},\qquad \theta\geq 0,
 \]
 is a Bernstein function where $\Gamma(u)$ denotes the usual Gamma function with parameter $u>0.$ In order to determine the triple $(\kappa, {\rm d}, \Upsilon)$ in (\ref{Bernstein}) associated with $\phi$ let us recall that the Beta function is related to the Gamma function by the following formula, for $a,b>0$  
$$B(a,b):=\int^1_{0}x^{b-1}(1-x)^{a-1} d x=\frac{\Gamma(a)\Gamma(b)}{\Gamma(a+b)}.$$
We thus have that $$\phi(\theta)=\frac{c\beta\theta}{\Gamma(\lambda)}B(\beta\theta+\nu,\lambda),\qquad \theta\geq 0.$$ Then making a change of variable in the expression for the Beta function we reach the identity 
\begin{equation}\label{remsubordinator}\begin{split}
\frac{\phi(\theta)}{\theta}=\frac{c}{\Gamma(\lambda)}\int^\infty_{0}e^{-\theta z}e^{-z\nu/\beta}\left(1-e^{-z/\beta}\right)^{\lambda-1} d z.
\end{split}
\end{equation}
This means that $\kappa={\rm d}=0$ and 
\[
\Upsilon(x,\infty)=c\frac{e^{-x\left(\nu+\lambda-1\right)/\beta}}{\Gamma(\lambda)}\left(e^{x/\beta}-1\right)^{\lambda-1},\qquad x>0.
\]
It is clear from the above expression that $\Upsilon$ has a density which is monotone decreasing.

In order to determine the potential measure associated to this subordinator observe the following elementary identity:
\begin{equation}\label{eq:ex1}\frac{\theta}{\phi(\theta)}=\frac{\Gamma(\nu+\beta\theta+\lambda)\Gamma(1-\lambda)}{c\beta\Gamma(\nu+\beta\theta+1)}\frac{\nu+\beta\theta}{\Gamma(1-\lambda)}.\end{equation} Therefore, we have that 
\begin{equation}\label{eq:kappa}
\begin{split}
\frac{\theta}{\phi(\theta)}
&=\frac{\nu+\beta\theta}{c\beta\Gamma(1-\lambda)}\int^1_{0} x^{\nu+\beta\theta-1}x^{\lambda}(1-x)^{-\lambda}  d x\\
&=\frac{\lambda}{c\beta\Gamma(1-\lambda)}\int^1_{0}\frac{ 1}{x^2}(1-x^{\nu+\beta\theta})\left(\frac{1}{x}-1\right)^{-\lambda-1}  d x\\
&=\frac{\lambda}{c\beta\Gamma(1-\lambda)}\int^\infty_{0} (1-e^{-(\nu+\beta\theta)z})\left(e^{z}-1\right)^{-\lambda-1}e^{z} d z\\
&=\frac{\lambda}{c\beta\Gamma(1-\lambda)}\int^\infty_{0} \left(1-e^{-\nu z}\right)\frac{e^{z}}{\left(e^{z}-1\right)^{\lambda+1}} d z\\
&\hspace{1cm} +\frac{\lambda}{c\beta\Gamma(1-\lambda)}\int^\infty_{0} \left(1-e^{-\beta\theta z}\right)\frac{e^{z(1-\nu)}}{\left(e^{z}-1\right)^{\lambda+1}} d z\\
&=\frac{\Gamma(\nu+\lambda)}{c\beta\Gamma(\nu)}+\frac{\lambda}{c\beta\Gamma(1-\lambda)}\int^\infty_{0} \left(1-e^{-\beta\theta z}\right)\frac{e^{z(1-\nu)}}{\left(e^{z}-1\right)^{\lambda+1}}  d z.\\
&=\frac{\Gamma(\nu+\lambda)}{c\beta\Gamma(\nu)}+\frac{\lambda}{c\beta^2\Gamma(1-\lambda)}\int^\infty_{0} \left(1-e^{-\theta x}\right)\frac{e^{x(1-\nu)/\beta}}{\left(e^{x/\beta}-1\right)^{\lambda+1}}  d x.
\end{split}
\end{equation}
This shows that $\phi$ is a special Bernstein function whose conjugate $\phi^*$ has triplet $(\kappa^*, {\rm d}^*, \Upsilon^*)$ (cf. (\ref{Bernstein*})) given by $\kappa^* = \Gamma(\nu+\lambda)/c\beta\Gamma(\nu)$, ${\rm d}^*=0$ and 
\[
\Upsilon^*(dx) = \frac{\lambda}{c\beta^2\Gamma(1-\lambda)}\frac{e^{x(1-\nu)/\beta}}{\left(e^{x/\beta}-1\right)^{\lambda+1}}  d x
\]
Note that $\Upsilon^*$  has a decreasing density.
Referring back to Theorem \ref{specialscale} and Corollary \ref{corrspecialscale} we may now say the following. 

There exists an oscillating spectrally negative L\'evy process with Laplace exponent
\[
 \psi(\theta) =  \frac{c\beta\theta^2\Gamma(\nu+\beta\theta)}{\Gamma(\nu+\beta\theta+\lambda)},\text{ for }
\theta\geq 0
\]
which has no Gaussian component and its L\'evy measure, $\Pi$, satisfies
\begin{eqnarray*}
 \Pi(-\infty,-x) &=& \frac{c(\nu+\lambda -1)}{\beta}\frac{e^{-x\left(\nu+\lambda-1\right)/\beta}}{\Gamma(\lambda)}\left(e^{x/\beta}-1\right)^{\lambda-1} \\
&&\hspace{2cm}+ \frac{c(\lambda-1)}{\beta}\frac{e^{-x\left(\nu+\lambda\right)/\beta}}{\Gamma(\lambda)}\left(e^{x/\beta}-1\right)^{\lambda-2} \text{ for }x>0
\end{eqnarray*}
and the associated scale function is given by 
\begin{equation}\label{gammascale}
W(x)= \frac{\Gamma(\nu+\lambda)}{c\beta\Gamma(\nu)}x+ \int_0^x\left\{ \int^\infty_{y}\frac{\lambda}{c\beta^2\Gamma(1-\lambda)}\frac{e^{z(1-\nu)/\beta}}{\left(e^{z/\beta}-1\right)^{\lambda+1}} dz \right\}dy
\end{equation}
for $x\geq 0.$

There exists a spectrally negative L\'evy process which drifts to $\infty$ or oscillates according to whether $\nu>0$ or $\nu=0,$ with Laplace exponent 
\begin{equation}\label{withoutcond}
\psi^*(\theta) = \frac{\theta\Gamma(\nu+\beta\theta+\lambda)}{c\beta\Gamma(\nu+\beta\theta)}\qquad \text{for}\ \theta\geq 0
\end{equation}
which has no Gaussian component and whose L\'evy measure, $\Pi^*$, satisfies 
\[
 \Pi^*(-\infty, -x) = \frac{\lambda}{c\beta^2\Gamma(1-\lambda)}\frac{e^{x(1-\nu)/\beta}}{\left(e^{x/\beta}-1\right)^{\lambda+1}}\qquad
\text{for}\ x>0
\]
and the associated scale function is given by 
\[
 W^*(x) =\int_0^x c\frac{e^{-z\left(\nu+\lambda-1\right)/\beta}}{\Gamma(\lambda)}\left(e^{z/\beta}-1\right)^{\lambda-1}dz.
\]

Note in the special case that $\nu=0$ and $\beta =c=1$ we have the two conjugate scale functions
\[
W(x) = \frac{1}{\Gamma(1-\lambda)}\int_0^x (e^y - 1)^{-\lambda}dy
\]
and
\[
W^*(x) = \frac{1}{\Gamma(\lambda)}\int_0^x (1- e^{-z})^{\lambda-1}dz.
\]

Another special case worthy of remark is the case where $\nu=1=\beta,$ $c=\Gamma(1+\lambda)$. The Laplace exponent $\psi^*$ takes the form $$
 \psi^*(\theta) = \frac{\Gamma(\theta +\alpha)}{\Gamma(\theta)\Gamma(\alpha)},\qquad \theta\geq 0,
$$ where $\alpha=1+\lambda\in(1,2).$ It was shown in Chaumont et al. \cite{CKP2007} that this is the Laplace exponent of a spectrally negative L\'evy process. The associated scale function, also identified in the latter paper, may simply be written 
\[
 W^*(x)  = (1-e^{-x})^{\alpha-1},\qquad x\geq 0.
\]
The previous calculations can be used to provide an example of the technique developed in Section \ref{genconst}. For $c,\nu>0,$ $\lambda\in]0,1[,$ let $\Psi$ be the function defined by means of $$\Psi(\theta)=(\theta-\nu)\frac{\Gamma(\theta+\lambda)}{c\Gamma(\theta)},\qquad\theta\geq 0.$$ It follows from the previous discussion, that the function $\Phi$ defined by $$\Phi(\theta)=\frac{\Gamma(\theta+\lambda)}{c\Gamma(\theta)}, \qquad \theta\geq 0,$$ is the Laplace exponent of a subordinator such that its L\'evy measure has a decreasing density and its potential measure has a decreasing and convex density.   According to our discussion in Section \ref{genconst} and the previous facts it follows that $\Psi(\theta)$ is the Laplace exponent of a spectrally negative L\'evy process that drifts to $-\infty,$ and its scale function is given by $$W(x)=\frac{c e^{\nu x}}{\Gamma(\lambda)}\int^x_{0}e^{-\nu z}(1-e^{-z})^{\lambda-1}dz,\qquad x\geq 0.$$ This is due to the fact that the L\'evy process with Laplace exponent $\Psi$ conditioned to drift to $\infty$ has Laplace exponent $$\Psi_{\nu}(\theta)=\Psi(\theta+\nu)=\frac{\theta\Gamma(\theta+\lambda+\nu)}{c\Gamma(\theta+\nu)},\qquad\theta\geq 0.$$ Moreover, the conjugate parent and ladder height process have Laplace exponent given by $$\Phi^*(\theta)=\frac{c\theta\Gamma(\nu+\theta)}{\Gamma(\nu+\lambda+\theta)},\qquad \Psi^*(\theta)=\frac{c\theta^2\Gamma(\nu+\theta)}{\Gamma(\nu+\lambda+\theta)},\qquad \theta\geq 0,$$ so that their characteristics were discussed before, and the corresponding scale function is described in equation (\ref{gammascale}), taking $\beta=1$. The Laplace exponent $\Psi$ can be taught as the one of the parent process with Laplace exponent $\psi^*,$ as in equation (\ref{withoutcond}), conditioned to drift to $-\infty.$ The particular case where $\nu=1=\beta,$ has been studied in \cite{CKP2007}. 

%\end{example}

\subsection*{Example 3}
%\begin{example}%[sum of stable subordinators]
\label{sumstable}
Let $0<\alpha\leq\beta\leq1$, $ a ,b>0$ and $\phi$ be the Bernstein function defined by 
$$
\phi(\theta)= a \theta^{\beta-\alpha}+b\theta^{\beta},\qquad \theta\geq 0.
$$ 
That is, in the case where $\alpha<\beta<1,$ $\phi$ is the Laplace exponent of a subordinator which is obtained as the sum of two independent stable subordinators one of parameter $\beta-\alpha$ and the other of parameter $\beta,$ respectively, so that the killing and drift term of $\phi$ are both equal to $0,$ and its L\'evy measure is given by $$\Upsilon(d x)=\left(\frac{ a (\beta-\alpha)}{\Gamma(1-\beta+\alpha)}x^{-(1+\beta-\alpha)}+\frac{b\beta}{\Gamma(1-\beta)}x^{-(1+\beta)}\right)d x,\quad x>0.$$ 

If $\alpha=\beta<1$ then $\phi$ is the Laplace exponent of a stable subordinator killed at rate $ a$;  when $\alpha<\beta=1,$ $\phi$ is the Laplace exponent of a stable subordinator with positive drift $b$; and finally in the case where $\alpha=1=\beta,$ $\phi$ is simply the Laplace exponent of a pure drift subordinator killed at rate $a$. The latter case  
will be excluded because it has been discussed in Example 1. In all cases the underlying L\'evy measure has a density which is completely monotone, and thus its potential density, or equivalently the density of the associated scale function $W$, is completely monotone. 

In the remainder of this example and subsequent examples we shall make heavy use of  the two parameter Mittag-Leffler function defined by 
\[
\mathrm{E}_{\alpha,\beta}(x) = \sum_{n\geq0}\frac{x^n}{\Gamma(n\alpha+\beta)},\qquad x\in\mathbb{R}.
\]
where $\alpha,\beta>0$.
The latter function can be identified via a pseudo-Laplace transform. Namely, for $\lambda\in\mathbb{R}$ and $\Re(\theta)>\lambda^{1/\alpha}-\gamma$,
\[
\int_0^\infty e^{-\theta x}e^{-\gamma x}x^{\beta-1}{\rm E}_{\alpha,\beta}(\lambda x^\alpha)dx=
\frac{(\theta+\gamma)^{\alpha-\beta}}{(\theta+\gamma)^\alpha-\lambda}.
\] 

The  associated scale function to $\phi$ can now be identified via 
\begin{equation} 
W'(x) = \frac{1}{b}x^{\beta-1}\mathrm{E}_{\alpha,\beta}\left(- a  x^{\alpha}/b\right),  \qquad x> 0,\label{MLreq}\end{equation} which is a completely monotone function because it is the product of the completely monotone functions $x^{\beta-1}$ and ${\rm E}_{\alpha,\beta}(-x^{\alpha}),$ and the later is completely monotone because it is the composition of the completely monotone function $t\mapsto {\rm E}_{\alpha,\beta}(-t)$ for $t\geq 0,$ see \cite{schneider96}, with the Bernstein function $x^\alpha.$ So, the function $$\psi(\theta)=\theta\phi(\theta)=a \theta^{\beta-\alpha+1} + b\theta^{\beta+1},  \qquad\theta\geq 0,$$ is the Laplace exponent of a spectrally negative L\'evy process. We shall elaborate on the features of the aforementioned parent process below according to three parameter regimes.

 In the case $\alpha<\beta<1,$ the parent process oscillates and is obtained by adding two independent spectrally negative stable processes with stability index $\beta+1$ and $1+\beta-\alpha,$ respectively. The scale function associated to it is given by 
$$
W(x)=\frac{1}{b}\int^x_{0}t^{\beta-1}\mathrm{E}_{\alpha,\beta}(- a  t^{\alpha}/b)d t,\qquad x\geq 0.
$$ 
The associated conjugates are given by $$\phi^*(\theta)=\frac{\theta}{a\theta^{\beta-\alpha}+b\theta^{\beta}},\quad \psi^*(\theta)=\frac{\theta^2}{ a \theta^{\beta-\alpha}+b\theta^{\beta}},\qquad \theta\geq 0,$$ and 
$$
W^*(x)=\frac{ a }{\Gamma(2-\beta+\alpha)}x^{1-\beta+\alpha}+\frac{b}{\Gamma(2-\beta)}x^{1-\beta},\qquad x\geq 0.$$ 
The subordinator with Laplace exponent $\phi^*$ has zero killing and drift terms and its L\'evy measure is obtained by taking the derivative of the expression in (\ref{MLreq}). By Theorem \ref{inverse} the spectrally negative L\'evy process with Laplace exponent $\psi^*,$ oscillates, has unbounded variation, has zero linear and Gaussian terms, and its L\'evy measure is obtained by derivating twice the expression in (\ref{MLreq}).

In the case, $\alpha=\beta<1,$ the Laplace exponent $\psi$ takes the form 
$$
\psi(\theta) = \theta \phi(\theta) =  a \theta + b\theta^{\beta+1}.
$$
The latter is the Laplace exponent of a oscillating spectrally negative $\alpha$-stable process with stability index $\alpha=(1+\beta),$ and positive drift with rate $ a $.
The scale function can be implicitly found in Furrer \cite{Fur1998} and it takes the form 
$$
W(x) = \frac{1}{ a }(1- \mathrm{E}_{\beta, 1}(- a  x^{\beta}/b)),\qquad x\geq 0.
$$ The  respective conjugates are given by  
\[
\phi^*(\theta)=\frac{\theta}{ a \theta+b\theta^{\beta}},\qquad \psi^*(\theta) = \frac{\theta^2}{ a  \theta+ b\theta^{\beta}},\qquad \theta\geq 0,
\] and  
\[
W^*(x) =  a  x + \frac{b}{\Gamma(2-\beta)} x^{1-\beta},\quad x\geq 0. 
\] The conjugate subordinator and spectrally negative L\'evy process, can be described using similar reasoning to that of the previous parameter regime and we omit the details. One may mention  here that by letting  $ a\downarrow 0 $   the Continuity Theorem for  Laplace transforms  tells us  that for the case  $\phi(\theta)=b\theta^\beta,$ the associated $\psi$ is the Laplace exponent of a spectrally negative stable process with stability parameter $1+\beta,$ and its scale function is given by $$W(x)=\frac{1}{b\Gamma(1+\beta)}x^{\beta},\qquad x\geq 0.$$ The associated conjugates are given by $$\phi^*(\theta)=b^{-1}\theta^{1-\beta},\qquad \psi^*(\theta)=b^{-1}\theta^{2-\beta},\qquad \theta\geq 0,$$ and $$W^*(x)=\frac{b}{\Gamma(2-\beta)} x^{1-\beta},\quad x\geq 0.$$ So that $\phi^*,$ respectively $\psi^*,$ corresponds to a stable subordinator of parameter $1-\beta,$ zero killing and drift terms; respectively, to a oscillating spectrally negative stable L\'evy process with stability index $2-\beta,$ and so its L\'evy measure is given by   
\[
\Pi^*(-\infty, -x) = \frac{\beta(1-\beta)}{b\Gamma(1+\beta)}x^{\beta-2},\qquad x\geq 0.
\]

Lastly, in the case $\alpha<\beta=1,$ the Laplace exponent of the parent process, $\psi$, is associated to the addition of a spectrally negative stable process with stability index $2-\alpha$ plus an independent continuous L\'evy process with no drift and Gaussian coefficient equal to $\sqrt{b}.$ Its associated scale function is given by $$W(x)=\frac{1}{b}\int^x_{0}\mathrm{E}_{\alpha,1}(- a  t^{\alpha}/b)d t,\qquad x\geq 0.$$ The respective conjugates are given by $$\phi^*(\theta)=\frac{\theta}{ a \theta^{1-\alpha}+b\theta},\quad \psi^*(\theta)=\frac{\theta^2}{ a \theta^{1-\alpha}+b\theta},\qquad \theta\geq 0,$$ and $$W^*(x)=b+\frac{ax^\alpha}{\Gamma(1+\alpha)},\qquad x\geq 0.$$ The characteristics of conjugate parent process can be determined as in the other two parameter regimes but what is different is that it is a process with bounded variation since $W^*(0)={\rm d}=b>0.$   

To complete this example, observe that the change of measure introduced in Lemma \ref{exptilting0} allows us to deal with the Bernstein function $$\phi(\theta)=k(\theta+m)^{\beta-\alpha}+b(\theta+m)^\beta,\quad \theta\geq 0,$$ where $m\geq 0$ is a fixed parameter. In this case we get that there exists a spectrally negative L\'evy process whose Laplace exponent is given by $$\psi(\theta)=k\theta(\theta+m)^{\beta-\alpha}+b\theta(\theta+m)^\beta,\quad \theta\geq 0,$$ and its associated scale function is given by $$W(x)=\frac{1}{b}\int^x_{0}e^{-mt}t^{\beta-1}\mathrm{E}_{\alpha,\beta}(- a  t^{\alpha}/b)d t,\qquad x\geq 0.$$ The respective conjugates can be obtained explicitly but we omit the details given that the expressions found are too involved.  

As in the Example 2, the degree of generality on which this example has been developed allows us to provide another example of the technique developed in Section \ref{genconst}. For, $m,a,b>0,$ $0<\alpha\leq \beta\leq 1,$  there exists a parent process drifting to $-\infty$ and with Laplace exponent 
$$\Psi(\theta)=(\theta-m)\left(a\theta^{\beta-\alpha}+b\theta^\beta\right),\qquad \theta\geq 0.$$ It follows from the previous calculations that the scale function associated to the parent process with Laplace exponent $\Psi$ is given by 
$$W(x)=\frac{e^{mx}}{b}\int^x_{0}e^{-mt}t^{\beta-1}\mathrm{E}_{\alpha,\beta}(- a  t^{\alpha}/b)d t,\qquad x\geq 0.$$

Finally observe that the function defined in (\ref{MLreq}) is a completely monotone function, which would have allowed us to present this example performing the construction indicated in the Remark 3 in Section \ref{completeclass}. 
%\end{example}

\subsection*{Example 4}
%\begin{example}
This example builds on the work of Boxma and Cohen \cite{BC1998} and the generalization thereof by Abate and Whitt \cite{AW1999}. In the latter, a parent process is considered whose Laplace exponent satisfies
\[
 \psi(\theta) = \theta - \frac{\lambda\theta}{(\mu+\sqrt{\theta})(1+\sqrt{\theta})}, \qquad \theta\geq 0.
\]
This corresponds to a process which has a linear unit drift minus a compound Poisson process of rate $\lambda$ with jumps whose distribution, $F$, has Laplace transform $1-\theta/(\mu+\sqrt{\theta})(1+\sqrt{\theta})$. Let $\eta(x) = e^x {\rm erfc}(\sqrt{x})$. The tail of the jump distribution satisfies 
\[
 F(x,\infty)=\left(\frac{1}{1-\mu}\right)(\eta(x)-\mu\eta(x\mu^2))
\]
which in the case $\mu=1$ should be interpreted in the limiting sense so that 
\[
 F(x,\infty)= (2x+1)\eta(x) - 2\sqrt{\frac{x}{\pi}}, \qquad x\geq 0.
\]
In both cases, the distribution also has mean $1/\mu$ and hence the mean of the L\'evy process is $1-\lambda/\mu$ which is assumed to be strictly positive (so that the process drifts to infinity). We may thus identify the characteristics of the Laplace exponent, $\phi$, of the descending ladder height process. Specifically $\kappa= 1/\mu$, ${\rm d}=0$ and $\Upsilon(dx) = \lambda F(x,\infty)dx$. The scale function associated with $\phi$ is given in \cite{AW1999} by
\[
W(x) =\frac{1}{1-\lambda/\mu}\left(1 - \frac{\lambda/\mu}{\nu_1 - \nu_2}(\nu_1 \eta(x\nu_2^2) - \nu_2\eta(x\nu_1^2))\right).
\]
where 
\[
\nu_{1,2} = \frac{1+\mu}{2}\pm\sqrt{\left(\frac{1+\mu}{2}\right)^2 - \left(1-\frac{\lambda}{\mu}\right)\mu}.
\]

Conveniently it is shown in \cite{AW1999} that $F(x,\infty)$  is completely monotone which makes $\phi$
 a complete Bernstein function. Moreover we automatically we get the existence of a conjugate complete scale function
 \[
 W^*(x) = \frac{1}{\mu}x +  \int_0^x \left\{\int_y^\infty\frac{(\eta(z)-\mu\eta(z\mu^2))}{1-\mu}dz\right\}dy
 \]
 (with the obvious interpretation when $\mu=1$). The associated parent process has Laplace exponent
 \[
 \psi^*(\theta)  =\frac{\theta^2(\mu+\sqrt{\theta})(1+\sqrt{\theta})}{(\mu-\lambda) +\theta +(1+\mu)\sqrt{\theta}}.
 \]
It is not difficult to show that $W(0+)=1$ and hence ${\rm d}^*=1$ showing that the conjugate parent process has a Gaussian component with coefficient $\sqrt{2}$. As usual, $\Pi^*(-\infty, -x)$  may be computed by considering $W'(x) - W'(+\infty)$ which happens, in this case, to be a rather cumbersome expression but, none the less, explicit.

\subsection*{Example 5}
%\begin{example}
Let $\alpha\in(0,1),$ $\kappa>0,$ ${\rm d}\geq 0,$ $c>0$ and $\phi$ be the Bernstein function defined by
 \begin{equation}
 \phi(\theta)=\kappa+{\rm d}\theta+c\theta^\alpha,\qquad \theta\geq 0.
\label{*}
 \end{equation}
 (Although the case ${\rm d}=0$ has been treated in Example 3 we include it here again because a different approach is proposed.) 
That is, $\phi$ is the Laplace exponent of a subordinator which is a $\alpha$-stable subordinator with drift ${\rm d}$ killed at rate $\kappa$. In this case 
$$\psi(\theta)=\kappa\theta+{\rm d}\theta^2+c\theta^{1+\alpha},\quad \theta\geq 0,$$ 
so that the parent L\'evy process associated to $\phi$ drifts to $\infty$ and is the sum of an independent Gaussian process with drift and a $(1+\alpha)$-stable process.  Given that the  L\'evy measure of $\phi$ has a density which is completely monotone it follows that its potential measure in $(0,\infty)$ has a density which is completely monotone,   which we will next describe. Let $X$ be an $\alpha$-stable subordinator, $\mathbf{e}_\kappa$ be an exponential random variable of parameter $\kappa,$ and assume that $X$ and $\mathbf{e}_\kappa$ are independent. Observe that the random variable $Z=\mathbf{e}_\kappa^{1/\alpha}X_{c}+{\rm d}\mathbf{e}_\kappa$ has a density and its Laplace transform has the form
$$\e(e^{-\theta Z})=\frac{\kappa}{\phi(\theta)},\qquad \theta\geq 0.$$ It follows that the scale function, $W$, associated to $\phi$  has density  given by $W'(x)=h_{Z}(x)/\kappa$ where $h_{Z}$ is the density of $Z$ and can be written  using that $X_{c}\stackrel{\text{law}}{=}c^{1/\alpha}X_{1}$ in terms of the density of $X_{1},$ say $p_{\alpha}$, or in terms of the Mittag-Leffler function, following if ${\rm d}>0$ or ${\rm d}=0,$ respectively, as follows 
$$h_{Z}(x)=\begin{cases}
\frac{\kappa}{c}\int^{xc/{\rm d}}_{0}e^{-\frac{\kappa}{c} s}p_{\alpha}\left(\frac{x-\frac{{\rm d}}{c}s}{s^{1/\alpha}}\right)\frac{d s}{s^{1/\alpha}},& \text{if} \ {\rm d}>0\\
(-1)\frac{d \mathrm{E}_{\alpha,1}\left(-\frac{\kappa}{c} x^{\alpha}\right)}{d x}=\frac{\kappa}{c}x^{\alpha-1}{\rm E}_{\alpha,\alpha}\left(-\frac{\kappa}{c}x^{\alpha}\right), & \text{if}\ {\rm d}=0,\end{cases}\qquad x>0;$$ owing to the fact that the Laplace transform of $X^{-\alpha}_{1},$ is given in terms of  the Mittag-Leffler function. An expression for the density $p_{\alpha}$ in series form can be found in equation (14.31) of Sato \cite{Sat1999}. It follows that the scale function associated to the spectrally negative L\'evy process with Laplace exponent $\psi$ is given by $\p(Z\leq x)/\kappa$. That is to say
\begin{eqnarray}
\hspace{-0.7cm}W(x)&=&\begin{cases}
\frac{c}{\kappa}\left(1-e^{-\kappa x/{\rm d}}\right)+\frac{c}{{\rm d}}e^{-\kappa x/{\rm d}}\int^x_{0}e^{\kappa u/{\rm d}}\p(X_{1}>u)du,& \text{if}\ {\rm d}>0\\
\frac{1}{\kappa}\left(1-\mathrm{E}_{\alpha,1}\left(-\frac{\kappa}{c} x^{\alpha}\right)\right),& \text{if} \ {\rm d}=0 
\end{cases} \label{use d=0}
\end{eqnarray}
for $x\geq 0$.
The conjugate Laplace exponent $\phi^*$ has zero drift and killing terms and its L\'evy measure is described by $\Upsilon^*(x,\infty)=W'(x),$ $x>0,$ the conjugate parent process oscillates and has Laplace exponent given by $$\psi^*(\theta)=\frac{\theta^2}{\kappa+{\rm d}\theta+c\theta^{\alpha}},\qquad \theta\geq 0.$$ Finally, the conjugate scale function is given by $$W^*(x)={\rm d}+\kappa x+\frac{c}{\Gamma(2-\alpha)}x^{1-\alpha},\qquad x\geq 0.$$ 

In the sequel % for $\kappa,{\rm d}\geq 0,$ $\alpha\in(0,1)$ we will denote by $\phi_{\kappa,{\rm d},\alpha}$ the Bernstein function $$\phi_{\kappa,{\rm d},\alpha}(\theta)=\kappa+{\rm d}\theta+\theta^{\alpha},\quad \theta\geq 0.$$  
fix $\gamma,c>0,$ ${\rm d},\kappa\geq 0$, $\alpha\in(0,1).$ We would like to determine the scale function associated to the Bernstein function 
\begin{equation}\label{stabletilted}\phi_{\gamma}(\theta)=\kappa+{\rm d}\theta+c(\theta+\gamma)^{\alpha}-c\gamma^{\alpha},\qquad \theta\geq 0.\end{equation}  To that end assume first that $\kappa> {\rm d}\gamma+c\gamma^{\alpha},$ and then observe that 
$$
\phi_{\gamma}(\theta)=\phi(\theta+\gamma), \quad \theta\geq 0,
$$ 
where $\phi$ has the form (\ref{*}) with killing term $\kappa-{\rm d}\gamma - c\gamma^\alpha$.
Note that $\phi_{\gamma}$ is the Laplace exponent of the subordinator with killing term $\kappa,$ drift term ${\rm d}$ and L\'evy measure 
$$
\Upsilon_{\gamma}(d x)=\frac{c\alpha e^{-\gamma x}}{\Gamma(1-\alpha)x^{1+\alpha}}d x,\qquad x>0.
$$ 
Note the density given above is the product of completely monotone functions and hence is itself completely monotone thus making $\phi_\gamma$ a complete Bernstein function.
The associated scale function, $W_\gamma$, has a density in $(0,\infty)$ which is given by $W'_{\gamma}(x) = e^{-\gamma x}W'(x)$. It follows from (\ref{use d=0}) and the fact that ${\rm d}_\gamma^*= \lim_{\theta\uparrow\infty}1/\phi_\gamma(\theta)=0,$ $\kappa^*=\lim_{\theta\to 0}\phi^*_{\gamma}(\theta)=\lim_{\theta\to 0}\frac{\theta}{\phi(\theta)}=0$,  
\begin{equation}
W_{\gamma}(x)=\begin{cases}
\frac{1}{c}\int^x_{0}e^{-\gamma y}\left(\int^{yc/{\rm d}}_{0}e^{-\frac{\kappa-{\rm d}\gamma-\gamma^\alpha}{c} s}p_{\alpha}\left(\frac{y-\frac{{\rm d}}{c}s}{s^{1/\alpha}}\right)\frac{d s}{s^{1/\alpha}}\right)dy,& \text{if} \ {\rm d}>0,\\
\frac{1}{c}\int_0^x e^{-\gamma y} y^{\alpha-1} {\rm E}_{\alpha,\alpha}\left(-\frac{(\kappa - c\gamma^\alpha)}{c}y^\alpha\right)dy,& \text{if} \ {\rm d}=0, 
\end{cases}
\label{COM}
\end{equation} for $x\geq 0.$ These also gives some insight about the form of the scale function for any value of $\kappa, \gamma,c>0$. Indeed, integrating directly in the case ${\rm d}>0,$ and  term by term in the case ${\rm d}=0,$ it is easily seen that  for any $\kappa,c>0,$ and $\gamma\geq 0$ the function 
\begin{equation*}
\begin{cases}
\frac{1}{c}e^{-\gamma y}\int^{yc/{\rm d}}_{0}e^{-\frac{\kappa-{\rm d}\gamma-\gamma^\alpha}{c} s}p_{\alpha}\left(\frac{y-\frac{{\rm d}}{c}s}{s^{1/\alpha}}\right)\frac{d s}{s^{1/\alpha}},& \text{if} \ {\rm d}>0,\\
\frac{1}{c} e^{-\gamma y} y^{\alpha-1} {\rm E}_{\alpha,\alpha}\left(-\frac{(\kappa - c\gamma^\alpha)}{c}y^\alpha\right),& \text{if} \ {\rm d}=0, 
\end{cases}\quad y\geq 0
\end{equation*}
has Laplace transform $1/\phi_{\gamma}(\theta),$ with $\phi_{\gamma}$ as defined in (\ref{stabletilted}). So, for any $\alpha\in(0,1),$ $\kappa,c,\gamma>0,$ ${\rm d}\geq 0,$ the scale function associated to $\phi_{\gamma}$ is given by (\ref{COM}). 

On account of the continuity theorem of Laplace transforms, the scale function $W$ is continuous in the L\'evy triple $(\kappa, {\rm d}, \Upsilon)$. For this reason, when ${\rm d}=0,$ taking further $\kappa\downarrow 0$, we recover from (\ref{COM}) 
\[
W_\gamma(x) =\frac{1}{c} \int_0^x e^{-\gamma y} y^{\alpha-1} {\rm E}_{\alpha,\alpha}( \gamma^\alpha y^\alpha)dy 
\] 
which is one of the scale functions found in \cite{HK2007}.% when one recalls that $\alpha {\rm E}_{\alpha,1}'(\lambda x^\alpha) = {\rm E}_{\alpha, \alpha}(\lambda x^\alpha)$ for all $\lambda \in \mathbb{R}$.

The conjugate scale function $W_\gamma^*$ is easily computed taking account of the L\'evy triple associated with $\phi_\gamma$ to be
\[
 W_\gamma^*(x) = {\rm d} + \kappa x + \frac{c\alpha }{\Gamma(1-\alpha)}\int_0^x \left\{\int_y^\infty   \frac{e^{-\gamma z}}{z^{1+\alpha}} dz \right\}dy
\]

Observe that taking $c=\lambda/\alpha$ in the definition of $\phi_{\gamma},$ in (\ref{stabletilted}), and making $\alpha$ tend to $0,$ we get that 
$$\lim_{\alpha\to 0}\left(\kappa+{\rm d}\theta+\frac{\lambda}{\alpha}\left((\theta+\gamma)^{\alpha}-\gamma^{\alpha}\right)\right)=k+{\rm d}\theta+\lambda\log\left(1+\frac{\theta}{\gamma}\right),\qquad \theta\geq 0.$$ So, by the continuity of Laplace transforms it follows that the scale function $W_{\gamma}$ converges as $\alpha\to 0$ to the scale function corresponding to the Bernstein function $k+{\rm d}\theta+\lambda\log\left(1+\frac{\theta}{\gamma}\right),$ which is studied in the case $\kappa=0={\rm d},$ $\gamma=1,$ in the forthcoming  Example 7.  

A different approach to this example would start with the Bernstein function in (\ref{use d=0}), which by the second remark following Corollary \ref{completescale} we know is the scale function associated to some spectrally negative L\'evy process, and then determine its associated parent process, ladder height process and conjugates. 
%\end{example}

\subsection*{Example 6}
%%%%%%%%%%%%%%%%%
%\begin{example}
 Another example belonging to a related family of Bernstein functions to those of the previous example starts again with calculations found in \cite{HK2007}. We take
\[
 \phi(\theta) = \kappa + \lambda\left(
1- \left(\frac{\gamma}{\gamma+\theta}\right)^\nu
\right),\qquad \theta\geq 0,
\]
where $\kappa,\lambda>0$. For this $\phi$ we also have ${\rm d}=0$ and 
\[
 \Upsilon(dx) = \frac{\lambda\gamma^\nu}{\Gamma(\nu)}x^{\nu-1}e^{-\gamma x}dx,\qquad x>0,
\]
where $\nu\in(0,1)$ and $\gamma>0$. Note that the assumption on $\nu$ ensures that $\Upsilon$ has a non-increasing density and hence $\phi$ may be used as a descending ladder height process. This Bernstein function is the Laplace exponent of a killed compound Poisson subordinator with Gamma distributed jumps. It can actually be seen as an extension to negative values for the parameter $\alpha$ in the definition of the Bernstein function considered in Example 5. It is also a complete Bernstein function on account of the fact that $\Upsilon$ has a completely monotone density. In \cite{HK2007} it was shown that 
\[
 W(x) =\frac{1}{\lambda+\kappa}  + \frac{\rho\gamma^\nu}{\lambda+\kappa} \int_0^x e^{-\gamma y}
  y^{\nu-1}{\rm E}_{\nu,\nu}(\rho\gamma^\nu y^\nu)dy
 % \frac{1}{\lambda+\kappa} e^{-\gamma x} {\rm E}_{\nu,1}(\rho\gamma^\nu x^\nu)
% + 
%\frac{\gamma}{\lambda+\kappa}\int_0^x e^{-\gamma y}{\rm E}_{\nu, 1}(\rho\gamma^\nu x^\nu)dy
\]
where $\rho= \lambda/(\lambda+\kappa)$.

Thanks to the fact that $\phi$ is a complete Bernstein function, its conjugate $\phi^*$ may also be used to construct a scale function. From the above description of $W$ one establishes  in a straightforward way that $\kappa^*=0$,
${\rm d}^* = 1/(\kappa+\lambda)$ and 
\[
\Upsilon^*(x,\infty) =\frac{\rho\gamma^\nu}{\lambda+\kappa}
e^{-\gamma x}  x^{\nu-1}{\rm E}_{\nu,\nu}(\rho\gamma^\nu x^\nu)
\]
and so 
\[
 W^*(x) = \kappa x + \frac{\lambda\gamma^\nu}{\Gamma(\nu)}\int_0^x \left\{\int_y^\infty  z^{\nu-1} e^{-\gamma z}dz\right\}dy.
\]

It follows that the respective conjugate parent processes have Laplace exponents
\[
 \psi(\theta) =\kappa\theta + \lambda\theta\left(
1- \left(\frac{\gamma}{\gamma+\theta}\right)^\nu
\right)\text{ and } \psi^*(\theta) = \frac{\theta^2}{\kappa + \lambda\left(
1- \left(\frac{\gamma}{\gamma+\theta}\right)^\nu
\right)}.
\]
Moreover, the parent process associated with $\phi$ drifts to $\infty$, has no Gaussian component and jump measure satisfying
\[
 \Pi(-\infty, -x) = \frac{\lambda\gamma^\nu}{\Gamma(\nu)}x^{\nu-1}e^{-\gamma x},\qquad \text{for}\ x>0.
\]
 The conjugate parent process oscillates, has a Gaussian coefficient $\sigma^*=\sqrt{2/(\lambda+\kappa)}$ and its jump measure satisfies
$$
\int^{\infty}_{x}\Pi^*(-\infty,-y)dy =\frac{\rho\gamma^\nu}{\lambda+\kappa} e^{-\gamma x}
  x^{\nu-1}{\rm E}_{\nu,\nu}(\rho\gamma^\nu x^\nu)
  ,\qquad x> 0.$$

Note that although it has been assumed that $\kappa>0$, since both scale functions are continuous in their parameters through their Laplace transform, it follows that one may simply take limits as $\kappa\downarrow 0$ to include $\kappa=0$ in the parameter ensemble. In the latter case, the parent process will oscillate as opposed to drifting to $\infty$.
%\[
% \Pi^*(-\infty, -x) = \frac{\rho^{1/\nu}\gamma}{\lambda+\kappa}e^{-\gamma x}\left(\gamma\frac{d}{dx}{\rm E}_{\nu, 1}(\rho\gamma^\nu x^\nu) -\frac{d^2}{dx^2}{\rm E}_{\nu, 1}(\rho\gamma^\nu x^\nu)\right) 
%\]

%\end{example}

\subsection*{Example 7}
%\begin{example}
\label{gamma-stable}
For $\alpha\in(0,1],$ $\lambda>0,$ now we take 
\[
\phi(\theta) = \lambda\log(1+\theta^\alpha),\qquad \theta\geq 0.
\]
In the case that $\alpha=1$, $\phi$ is the Laplace exponent of a Gamma subordinator. In the case that $\alpha\in(0,1)$ this is the Laplace exponent of an $\alpha$-subordinator, subordinated by a Gamma subordinator. A subordinator with $\phi$ as Laplace exponent is usually called Linnik subordinator.   
For this Laplace exponent one may show that $\kappa  = {\rm d}=0$ and 
\[
\Upsilon(dx) = \alpha\lambda \frac{{\rm E}_{\alpha, 1}(-x)}{x}dx,\qquad x>0.
\]
Note that $\phi$ is a complete Bernstein function as soon as one recalls that ${\rm E}_{\alpha,1}(-x)$ (cf. \cite{Pol1948}) and $1/x$ are completely monotone and that the product of completely monotone functions is completely monotone.
The scale function associated with $\phi$ is unknown. None the less, from the discussion in Section \ref{completeclass} we see that the conjugate $\phi^*$ has a L\'evy measure $\Upsilon^*$ which is absolutely continuous with non-increasing density. Hence it follows from Corollary \ref{corrspecialscale} that there exists a spectrally negative L\'evy process with Laplace exponent 
\[
\psi^*(\theta) = \frac{\theta^2}{\lambda\log(1+\theta^\alpha)},\qquad \theta\geq 0,
\]
whose associated scale function is
\[
W^*(x) =\alpha\lambda \int_0^x \left\{\int_y^\infty  \frac{{\rm E}_{\alpha, 1}(-z)}{z}dz \right\}dy,\qquad x\geq 0.
\]

It was shown in Hubalek and Kyprianou \cite{HK2007} that for the case $\alpha=1$ 
\[
W(x) = \lambda\int_0^x e^{-y}\left\{\int_0^\infty \frac{y^{t-1}}{\Gamma(t)}dt \right\}dy,\qquad x\geq 0.
\]
It follows from the latter that ${\rm d}^*=0$. Since $\phi^*(\theta)= \theta/\lambda\log(1+\theta)$, an easy limit as $\theta\downarrow 0$ shows that $\kappa^* = 1$ and hence  from the expression given for $W$
\[
\Upsilon^*(x,\infty) = \lambda e^{-x}\int_0^\infty \frac{x^{t-1}}{\Gamma(t)}dt-1,\qquad x>0.
\]
It now follows that the conjugate parent process drifts to $\infty$, has no Gaussian component and, formally, Theorem \ref{inverse} tells us that $\Pi^*$ satisfies
\[
\int^\infty_{x}\Pi^*(-\infty, -y)dy = \lambda e^{-x}\int_0^\infty \frac{x^{t-1}}{\Gamma(t)}dt,\qquad x>0.
\]
%\end{example}

\subsection*{Example 8}
%\begin{example}
This example is built again on discussion found in \cite{SV2007}. Suppose we take 
\[
\phi(\theta) = \log \{(1+\theta) +  \sqrt{(1+\theta)^2 -1}\},\qquad \theta\geq 0.
\]
This is a Bernstein function which has $\kappa={\rm d}=0$ and 
\[
\Upsilon(dx) = \frac{e^{-x}}{x}I_0(x)dx,\qquad x>0,
\] 
where $I_\nu$ is the modified Bessel function of  index $\nu$. It is known that the density of $\Upsilon$ is the Laplace transform of 
\[
\gamma(dy) = \left\{\int_0^y  \frac{1}{\pi}(2z - z^2)^{-1/2}\mathbf{1}_{\{z\in(0,2)\}} dz\right\}dy,\qquad y>0,
\]
and hence $\phi$ is a complete Bernstein function.
Following the reasoning in the previous example, we have the existence of an oscillating spectrally negative L\'evy process with Laplace exponent 
\[
\psi^*(\theta) = \frac{\theta^2}{ \log \{(1+\theta) +  \sqrt{(1+\theta)^2 -1}\}}
\]
for $\theta\geq 0,$
whose scale function is given by 
\[
W^*(x) =  \int_0^x \left\{\int_y^\infty  \frac{e^{-z}}{z}I_0(z)dz \right\}dy,\qquad x\geq 0.
\]

In \cite{SV2007} one also finds that the transition density of the subordinator, $H$, associated with $\phi$ is known and specifically is given by 
$$
\mathbb{P}(H_t \in dx) = \frac{t}{x}I_t(x)e^{-x}
$$
and hence from (\ref{resolvent2}) it follows that
\[
W(x) = \int_0^x \frac{e^{-y}}{y}\left\{\int_0^\infty t I_t (y)dt\right\} dy.
\]
Clearly we may now say that ${\rm d}^*=0$ and moreover from the fact that $\phi^*(\theta) = \theta/\phi(\theta)$, taking limits as $\theta\downarrow 0$ tells us also that $\kappa^*=0$. Moreover the expression for $W$ given above also gives us
\[
\Upsilon^*(x,\infty) =  \frac{e^{-x}}{x}\int_0^\infty t I_t (x)dt,\qquad x\geq 0.
\]
We see then that the conjugate parent process oscillates, has no Gaussian component and, formally, Theorem \ref{inverse} tells us that
\[
\int_x^\infty\Pi^*(-\infty, -y)dy = \frac{e^{-x}}{x}\int_0^\infty t I_t (x)dt,\qquad x>0.
\]

\section*{Acknowledgments}
We are grateful to Reinming Song and Zoran Vondra\v cek for bringing to our attention several bibliographical references and for insightful comments that allowed us to improve the presentation of this work.    
This work was carried out over various periods for which support was provided by the EPSRC grants
EP/D045460/1 and EP/C500229/1, and the CONCyTEG grant 06-02-K117-35.

\begin{flushleft}
Department of Mathematical Sciences\\ 
The University of Bath \\
Claverton Down \\
Bath BA2 7AY \\
UK. \\
email: a.kyprianou@bath.ac.uk
\end{flushleft}

\begin{flushleft}
Centro de Investigación en Matemáticas A.C. \\
Calle Jalisco s/n\\
Col. Valenciana\\
C.P. 36240 Guanajuato, Gto.\\
Mexico\\
email: rivero@cimat.mx\\
\end{flushleft}
\end{document}